\documentclass[12pt]{article}

\usepackage[english]{babel} 
\usepackage{graphicx}
\usepackage{fancyhdr}

\newtheorem{theorem}{Theorem}[section]
\newtheorem{proposition}{Proposition}[section]
\newtheorem{lemma}[theorem]{Lemma}
\newtheorem{corollary}{Corollary}[section]
\newtheorem{example}[theorem]{Example}
\newtheorem{remark}[theorem]{Remark}
\newtheorem{proof}{Proof}

\usepackage{amsmath}
\usepackage{amssymb}
\usepackage{pstricks}

\def\proof{{\it Proof.}\ }
\def\QED{\hfill $\square$}

\date{}

\title{Factorization formulas for Macdonald polynomials} 

\author{Fran\c{c}ois Descouens \scriptsize (1) \normalsize and Hideaki Morita \scriptsize (2) \normalsize \\
  \scriptsize (1) Institut Gaspard Monge, Universit\'{e} de Marne-la-Vall\'{e}e 77454 Marne-la-Vall\'{e}e, France  \\ 
                 \scriptsize               francois.descouens@univ-mlv.fr\\
      \scriptsize         (2) Oyama National College of Technology, Oyama, Tochigi, 323-0806, Japan \\ 
                 \scriptsize           morita@oyama-ct.ac.jp  \normalsize
       }

\begin{document}

\maketitle

\begin{abstract}
  The aim of this note is to give some factorisation formulas for
  different versions of the Macdonald polynomials when the parameter
  $t$ is specialized at roots of unity, generalizing those given in
  \cite{LLT1} for Hall-Littlewood functions.
\end{abstract}

\section{Introduction}
In \cite{LLT1}, Lascoux, Leclerc and Thibon give some factorisation
formulas for Hall-Littlewood functions when the parameter $q$ is
specialized at roots of unity. They also give formulas in terms of
cyclic characters of the symmetric group. In this article, we give a
generalization of these specializations for different versions of the
Macdonald polynomials and we obtain similar formulas in terms of
plethysms and cyclic characters. In the last section, we give
congruence formulas for $(q,t)$-Kostka polynomials using Schur
functions in the alphabet of the powers of the
parameter $t$. We will mainly follow the notations of \cite{macd}.  \\
\\\noindent {\bf Acknowlegdements:} All computations on Macdonald
polynomials have been done using the {\tt MuPAD} package {\tt
  MuPAD-Combinat} (see \cite{HT} for more details on the project and
the website {\it http://mupad-combinat.sourceforge.net/}).

\section{Preliminaries}
For a partition $\lambda=(\lambda_1,\ldots,\lambda_n)$, $\lambda_1\ge
\ldots \ge \lambda_1$, we write $l(\lambda)$ its length, $\vert \lambda
\vert$ its weight, $m_i(\lambda)$ the multiplicity of the part of
length $i$ and $\lambda^{'}$ its conjugate partition.  Let $q$ and $t$
be two indeterminates and $F={\mathbb{Q}}(q, t)$.  Let $\Lambda_F$ be
the ring of symmetric functions over the field $F$.  Let denote by
$\langle\cdot,\cdot\rangle_{q,t}$ the inner product on $\Lambda_F$
defined on the powersums by
$$
\langle\  p_{\lambda}\ ,\ p_{\mu} \ \rangle_{q,t} =
\delta_{\lambda\mu} z_{\lambda}(q,t),
$$
where
$$
z_{\lambda}(q,t) = \prod_{i\ge 1}(m_i)!\ i^{m_i(\lambda)}\ 
\prod_{i=1}^{l(\lambda)} \frac{1-q^{\lambda_i}}{1-t^{\lambda_i}}.
$$
The special case $ \langle\cdot,\cdot\rangle :=
\langle\cdot,\cdot\rangle  {}_{q=0\atop t=0} $ is the usual inner
product.  

\noindent 
Let $\lbrace P_{\lambda}(x; q,t)\rbrace_\lambda$ be the family of
Macdonald polynomials obtained by orthogonalization of the Schur
functions basis with respect to the inner product
$\langle\cdot,\cdot\rangle_{q,t}$. Let us define a normalization of these
functions by
$$
Q_{\lambda}(x; q, t) = \frac{1}{\langle P_{\lambda}(x; q, t),
  P_{\lambda}(x; q, t)\rangle_{q,t}} P_{\lambda}(x; q,t).
$$
It is clear from the previous definitions that the families
$\lbrace P_{\lambda}(x;q,t)\rbrace _\lambda $ and $\lbrace
Q_{\mu}(x;q,t)\rbrace_\mu$ are dual to each other with respect to the
inner product $\langle\cdot,\cdot\rangle_{q,t}$ (c.f. \cite{macd},
Chap. I, section 4 and Chap. VI, formula (2.7)]).


\begin{proposition}\label{CauchyPQ}{\em [\cite{macd},VI, (4.13)] \em}
  Let $x=(x_1,x_2,\dots )$ and $y=(y_1,y_2,\dots )$ be two alphabets.
  The Macdonald polynomials $\lbrace P_\lambda(x;q,t) \rbrace_\lambda
  $ and $\lbrace Q_\lambda(x;q,t) \rbrace_\lambda $ satisfy the Cauchy
  formula
  \begin{eqnarray}
  \sum_{\lambda}P_{\lambda}(x; q, t)Q_{\lambda}(y; q, t) =
  \prod_{i,j} \frac {(tx_iy_j;q)_{\infty}} {(x_iy_j;q)_{\infty}}\ ,
  \end{eqnarray}
  where $(a;q)_{\infty}$ is the infinite product $ \prod_{r\geq
    0}(1-aq^r) $.
\end{proposition}


\noindent We consider the following algebra homomorphism
$$
\begin{array}{lll}
        {}' : \Lambda_F & \longrightarrow & \Lambda_F \\ 
                   f(x) & \longmapsto     & f'(x)=f\left(\frac{1-q}{1-t}x\right).
\end{array}
$$
The images of the powersums $(p_k)_k\ge 1$ by these morphisms are
$$ 
p_k^{'}(x)=\frac{1-q^k}{1-t^k}p_k(x).
 $$
 
\noindent
Let us consider the following modified version of the Macdonald polynomial
$$
Q'_{\mu}(x; q, t) = Q_{\mu} \left( \frac{1-q}{1-t}x; q, t \right).
$$
We can see that the families $\{Q'_{\mu}(x;q,t)\}_\mu$ and
$\{P_{\lambda}(x;q,t)\}_\lambda$ are dual to each other with respect to the
usual inner product.


\begin{proposition}\label{CauchyPQP}
  Let $x=(x_1,x_2,\dots )$ and $y=(y_1,y_2,\dots )$ be two alphabets.
  The Macdonald polynomials $\lbrace P_\lambda(x;q,t) \rbrace_\lambda $ and
  $\lbrace Q^{'}_\lambda(x;q,t) \rbrace_\lambda$ satisfy the following Cauchy
  formula
  \begin{eqnarray}
  \sum_{\lambda} P_{\lambda}(x; q, t) Q'_{\lambda}(y; q, t) =
  \prod_{i, j} \frac{1}{1-x_iy_j}.
  \end{eqnarray}
\end{proposition}

\noindent
\proof
By Proposition \ref{CauchyPQ}, we have
$$
\sum_{\lambda}P_{\lambda}(x; q, t)Q_{\lambda}(y; q, t) =
\prod_{i,j} \prod_{r\geq 0} \frac {1-tx_iy_jq^r} {1-x_iy_jq^r}.
 $$
 Since the map $x\mapsto x/(1-t)$ corresponds to the transformation
 of the alphabet $\{x_1,x_2,\ldots\}$ into the alphabet $\{x_it^j,\ i\ge 1, j\ge
 0\}$, a straightforward computation shows that
 $$
 \sum_{\lambda} P_{\lambda}(x; q, t)
 Q_{\lambda}\left(\frac{y}{1-t}; q, t\right) = \prod_{r\geq
   0}\prod_{i,j} \frac{1}{1-x_iy_jq^r}.
 $$
 This means that the families $\{P_{\lambda}(x; q, t)\}_\lambda$ on
 the alphabet $x$ and $\left \{Q_{\mu}\left(\frac{y}{1-q}; q, t\right)\right \}_\mu$ on the alphabet
 $y/(1-q)$ are dual to each other with respect to the usual inner
 product.  Since the transformation of alphabets $y\mapsto y/(1-q)$ is
 invertible and the inverse map is given by $y\mapsto (1-q)y$, it
 follows that
 $$
 \sum_{\lambda} P_{\lambda}(x; q, t)
 Q_{\lambda}\left(\frac{1-q}{1-t}y; q, t\right) =P_{\lambda}(x; q, t)
 Q'_{\lambda}(y; q, t)= \prod_{i,j}\frac{1}{1-x_iy_j}.
$$
\QED

\noindent We recall some definitions of combinatorial quantities associated to a
cell $s=(i,j)$ of a given partition. The arm length $a(s)$,
arm-colength $a^{'}(s)$, leg length $l(s)$ and leg-colength $l^{'}(s)$
are respectively the number of cells at the east, at the west, at the
south and at the north of the cell $s$ (c.f. \cite{macd}, Chap. VI, formula (6.14)), i.e
\begin{eqnarray*}
  a(s)= \lambda_i-j~    & , & ~ a^{'}(s)=j-1, \\
  l(s)= \lambda^{'}_j-i & , & ~ l^{'}(s)=i-1 .
\end{eqnarray*}
We also define the quantity
$$  n(\lambda)=\sum_{i}(i-1)\lambda_i.$$
\noindent Let
$ J_{\mu}(x; q, t) $ be the symmetric function with two parameters
defined by
\begin{eqnarray}\label{IntegralFormDef}
J_{\mu}(x; q, t) = c_{\mu}(q, t)P_{\mu}(x; q,t) = c'_{\mu}(q,
t)Q_{\mu}(x; q, t)\ , 
\end{eqnarray}
with
$$
c_{\mu}(q, t)=\prod_{s\in\mu}(1-q^{a(s)}t^{l(s)+1})\quad \text{and}
\quad c'_{\mu}(q, t)=\prod_{s\in\mu}(1-q^{a(s)+1}t^{l(s)}).
$$
The symmetric function $J_{\mu}(x; q, t)$ is called the integral
form of $P_{\mu}(x; q,t)$ or $Q_{\mu}(x; q, t)$ (c.f. \cite{macd}, Chap. VI, section
8).  Using this integral form, we can define an other modified
version of the Macdonald polynomial and the $(q,t)$-Kostka polynomials
$K_{\lambda, \mu}(q,t)$ by
\begin{eqnarray}
J_{\mu} \left( \frac{x}{1-t}; q, t
\right)= \sum_{\lambda} K_{\lambda,\mu}(q, t) s_\lambda.
\end{eqnarray}
In \cite{HHL}, Haglund, Haiman and Loehr consider a modified version
of $J_{\mu} \left( \frac{x}{1-t}; q, t
\right)$ and introduce other $(q,t)$-Kostka
polynomials $\tilde{K}_{\lambda, \mu}(q,t)$ by defining the functions
\begin{eqnarray}
\tilde{H}_{\mu}(x; q,t) = t^{n(\mu)}J_{\mu} \left( \frac{x}{1-t^{-1}}; q, t^{-1}
\right)=
\sum_{\lambda} \tilde{K}_{\lambda,\mu}(q, t) s_\lambda\ .
\end{eqnarray}
They give a combinatorial interpretation of this modified version
expanded on the monomials basis by defining two statistics (major index and
inversions) on arbitrary fillings by integers of $\mu$.

\begin{remark}
  {\em Let $\mu$ and $\rho$ be two partitions of the same weight. We have
    $$
    X^{\mu}_{\rho}(q, t) = \langle \ \tilde{H}_{\mu}(q, t)\ ,\ 
    p_{\rho}(x)\ \rangle,
    $$
    where $X^{\mu}_{\rho}(q, t)$ is the Green polynomial with two
    variables, defined by
    $$
    X^{\mu}_{\rho}(q, t) =
    \sum_{\lambda}\chi^{\lambda}_{\rho}\tilde{K}_{\lambda\mu}(q, t).
    $$
    Here $\chi^{\lambda}_{\rho}$ is the value of the irreducible
    character of the symmetric group corresponding to the partition
    $\lambda$ on the conjugacy class indexed by $\rho$. For related topics, see \cite{Morita1, Morita2}. \em}
\end{remark}


\section{Plethystic formula}
\noindent In this section, we prove a plethystic formula for
Macdonald polynomials indexed by rectangular partitions when the second
parameter $t$ is specialized at primitive roots of unity.
\begin{proposition}\label{Alphabett}{\em [\cite{macd},VI, (6.11')] \em}
  Let $l$ be a positive integer and $\lambda$ be a partition such that
  $l(\lambda)\le l$.  The Macdonald polynomials $P_\lambda(x;t,q)$ on
  the alphabet $\{x_i=t^i,\ 0\le i \le l-1, \ \text{ and } x_i=0,\ \forall
  i\ge l\}$, can be written
 \begin{eqnarray}\label{PAlphabet}
        P_{\lambda}(1, t, \dots ,t^{l-1}; q, t)
  =
  t^{n(\lambda)}
  \prod_{s\in\lambda}
  \frac{1-q^{a'(s)}t^{l-l'(s)}}{1-q^{a(s)}t^{l(s)+1}}.
 \end{eqnarray}
\end{proposition}

\begin{corollary}\label{Alphabetroot}
  Let $l$ be a positive integer and $\lambda$ a partition such that
  $l(\lambda) \le l$. For $\zeta$ a primitive $l$-th root of unity,
  the Macdonald polynomials $P_\lambda(x;q,t)$ satisfy the
  specialization
  \begin{eqnarray}
  P_{\lambda}(1, \zeta, \zeta^2, \dots ,\zeta^{l-1}; q, \zeta) =
  \left\{
        \begin{array}{cc}
        (-1)^{(l-1)r}
        &\mbox{if $\lambda=(r^l)$ for some $r\geq 0$}, \\
     0 &\mbox{otherwise}.
    \end{array}
   \right.
 \end{eqnarray}
\end{corollary}

\noindent
\proof Supplying zeros at the end of $\lambda$, we consider the
  partition $\lambda$ as a sequence of length exactly equal to $l$.
  The multiplicity of $0$ in $\lambda$ is $m_0=l-l(\lambda)$. We will
  denote by $\varphi_r(t)$ the polynomial $$
  \varphi_r(t)=(1-t)(1-t^2)\ldots (1-t^r).$$

\noindent Let
$$
        f(t)
        =
        \frac
         {(1-t^{l})(1-t^{l-1})\cdots(1-t^{l-l(\lambda)})
        (1-t^{l-l(\lambda)-1})\cdots(1-t^2)(1-t)}
         {\varphi_{m_0}(t)\varphi_{m_1}(t)\varphi_{m_2}(t)\cdots\cdots}
         $$
         be the product of factors of the form $1-q^0t^{\alpha}$
         for some $\alpha>0$ in the formula (\ref{PAlphabet}).  If we suppose that
         $f(\zeta)\neq 0$, the factor $1-t^l$ should be contained in
         one of $\varphi_{m_i}(t)$. This means that there exists
         $i\geq 0$ such that $m_i\geq l$.  Since we consider $\lambda$
         as a sequence of length exactly $l$, this implies the
         condition $m_r=l$ for some $r\geq 0$.  Thus, if
         $P_{\lambda}(1,\zeta,\zeta^2,\dots, \zeta^{l-1}; q,
         \zeta)\neq 0$, the shape of $\lambda$ should be $(r^l)$.  \\ 
         \\ Suppose now that $\lambda=(r^l)$.  By Proposition \ref{Alphabett}, it follows
         that
\begin{eqnarray*}
        P_{\lambda}(1, \zeta, \zeta^2, \dots ,\zeta^{l-1}; q, \zeta)
                &=&
                        \zeta^{n(\lambda)}
                        \prod_{s\in \lambda}
                        \frac
                                {1-q^{a'(s)}\zeta^{l-l'(s)}}
                                {1-q^{a(s)}\zeta^{1+l(s)}}\\
                &=&
                        \zeta^{n(\lambda)}
                        \prod_{(i,j)\in \lambda}
                        \frac
                                {1-q^{j-1}\zeta^{l-(i-1)}}
                                {1-q^{r-j}t^{l-i+1}}\\
                &=&
        \zeta^{n(\lambda)}
                        \prod_{i=1}^{l}
   \prod_{j=1}^r
                        \frac
                                {1-q^{j-1}\zeta^{l-i+1}}
                                {1-q^{r-j}t^{l-i+1}}\ .
\end{eqnarray*}
For each $i$, 
it is easy to see that 
$$
        \prod_{j=1}^r
                        \frac
                                {1-q^{j-1}\zeta^{l-i+1}}
                                {1-q^{r-j}t^{l-i+1}}
 =1.
$$ 
Hence, we obtain 
$$
        P_{\lambda}(1, \zeta, \zeta^2, \dots ,\zeta^{l-1}; q, \zeta)
 =
 \zeta^{n(\lambda)}, 
$$
and it follows immediately from the definition of $n(\lambda)$ 
that 
$$
        \zeta^{n(\lambda)}
 =
 \zeta^{l(l-1))r/2}
 =
 (-1)^{(l-1)r}.
$$
\QED

\begin{theorem}\label{QProot}
  Let $l$ and $r$ be two positive integers and $\zeta$ a primitive
  $l$-th root of unity. The Macdonald polynomials $Q^{'}_{r^l}(x;q,t)$
  satisfy the following specialization formula at $t=\zeta$
        \begin{eqnarray}
                Q'_{(r^l)}(x; q, \zeta) &=& (-1)^{(l-1)r} (p_l\circ
                h_r)(x)\ .
        \end{eqnarray}
\end{theorem}

\noindent
\proof
  Recall that
  $$
  \sum_{\lambda} P_{\lambda}(x; q,t) Q'_{\lambda}(y; q,t) =
  \prod_{i,j} \frac{1}{1-x_iy_j}.
  $$
  If we let $x_i=\zeta^{i-1}$, for $i=1,2,\dots ,l$, and $x_i=0$, for
  $i>l$, and $t=\zeta$, we obtain
  \begin{eqnarray}\label{CauchyRoot}
  \sum_{\lambda} P_{\lambda}(1, \zeta, \zeta^2, \dots, \zeta^{l-1};
  q,\zeta) Q'_{\lambda}(y; q,\zeta) = \prod_{j\geq 1} \prod_{i=1}^l
  \frac{1}{1-\zeta^{i-1}y_j}.
  \end{eqnarray}
  By Corollary \ref{Alphabetroot}, the left hand side of (\ref{CauchyRoot}) is
  equal to
  $$
  \sum_{r\geq 0} (-1)^{(r-1)l}Q'_{(r^l)}(y; q, \zeta).
  $$
  Since $\prod_{i=1}^l(1-\zeta^{i-1}t)=1-t^l$, the right hand side
  of (\ref{CauchyRoot}) coincides with
  $
  \sum_{r\geq o}h_r(y^l)
  $, 
  where $y^l$ denotes the alphabet $(y_1^l, y_2^l,
  \cdots)$.  Comparing the degrees, we can conclude that
  $$
  Q'_{(r^l)}(y; q, \zeta) = (-1)^{(l-1)r}h_r(y^l) =
  (-1)^{(l-1)r}(p_l\circ h_r)(y).
  $$
\QED
\begin{example}{\em For $\lambda=(222)$ and $l=3$, we can compute the specialization
\begin{eqnarray*} 
 Q_{(222)}^{'}(x;q,e^{\frac{2i\pi}{3}})&=& - s_{321} + s_{33} + s_{411} - s_{51} + s_{6} + s_{222}\\
 &=&p_3\circ h_2(x) .
\end{eqnarray*} \em}
\end{example}

\noindent In order to give a similar formula for the modified versions 
of the integral form of the Macdonald polynomials, we give a formula
for the specialization of the constant $c'_{r^l}(t,q)$ at $t$ a
primitive $l$-th root of unity.
\begin{lemma}\label{Constantroot} 
  Let $l$ and $r$ be two positive integers and $\zeta$ a $l$-th
  primitive root of unity. The normalization constant
  $c_\lambda^{'}(q, t)$ satisfies the following specialization at
  $t=\zeta$
\begin{eqnarray} 
c_{r^l}^{'}(q, \zeta)= \prod_{i=1}^{r}(q^{il}-1).
\end{eqnarray}

\end{lemma}
\noindent
\proof Recall the definition of the normalization constant
\begin{eqnarray*}
c'_{r^l}(q, t)=\prod_{s\in\mu}(1-q^{a(s)+1}t^{l(s)})=\prod_{i=1}^{r}\prod_{j=1}^{l}(1-q^{r-i+1}t^j)=
\prod_{i=1}^{r}\prod_{j=1}^{l}(1-q^{i}t^j).
\end{eqnarray*}
Specializing $t$ at $\zeta$ a $l$-th primitive root of unity, we obtain
\begin{eqnarray*}
c'_{r^l}(q, \zeta)=\prod_{i=1}^{r}\prod_{j=1}^{l}(1-q^{r-i+1}t^j)=\prod_{i=1}^{r}(q^{il}-1).      
\end{eqnarray*}
\QED
\begin{corollary}
  With the same notations as in Theorem \ref{QProot}, the modified
  integral form of the Macdonald polynomials $\tilde{H}_\mu(x;q,t)$
  satisfy a similar formula at $t=\zeta$,
  \begin{eqnarray}
 \tilde{H}_{(r^l)}(x;q,\zeta)&=&
 \prod_{i=1}^{r}(q^{il}-1)\ p_l\circ h_r\left (\frac{x}{1-q}\right ).
\end{eqnarray}
\end{corollary}

\noindent
\proof  Using the definition \ref{IntegralFormDef} of the integral form of the
Macdonald polynomials
\begin{eqnarray*}
J_{(r^l)}\left (\frac{x}{1-t};q,t\right )
=& c_{{r^l}}^{'}(q,t)Q_{(r^l)}\left (\frac{x}{1-t};q,t\right ) \\
=& c_{{r^l}}^{'}(q,t)Q^{'}_{{(r^l)}}\left (\frac{x}{1-q};q,t\right )\ .
\end{eqnarray*}
This expression can be rewritten in terms of plethysms by the powersum
$p_1$, consequently
$$J_{(r^l)}\left (\frac{x}{1-t};q,t\right )=c_{{r^l}}^{'}(q,t)\left (
  Q^{'}_{(r^l)}(\ .\ ;q,t)\circ \frac{1}{1-q}p_1\right )(x)\ . $$
By
specializing in this egality, $t$ at a primitive $l$-th root of unity
$\zeta$, using Theorem \ref{QProot}, we obtain
\begin{eqnarray*}
J_{(r^l)}\left (\frac{x}{1-\zeta};q,\zeta\right )&=& c_{{r^l}}^{'}(q,\zeta)\left (
  Q^{'}_{(r^l)}(\ .\ ;q,\zeta )\circ \frac{1}{1-q}p_1\right )(x)\\
&=& c_{{r^l}}^{'}(q,\zeta)(-1)^{r(l-1)}\left (\left ( p_l \circ h_r\right ) \circ \frac{1}{1-q}p_1\right )(x).
\end{eqnarray*}
As plethysm is associative, we can write
\begin{eqnarray*}
J_{(r^l)}\left (\frac{x}{1-\zeta};q,\zeta\right )
&=&c_{{r^l}}^{'}(q,\zeta)(-1)^{r(l-1)}p_l \circ \left ( h_r \circ \frac{1}{1-q}p_1\right )(x) \\
&=& c_{{r^l}}^{'}(q,\zeta)(-1)^{r(l-1)}p_l \circ h_r\left (\frac{x}{1-q} \right ). 
\end{eqnarray*}
Using the formula of Lemma \ref{Constantroot} and $
\zeta^{n(\lambda)}=(-1)^{(l-1)r}.$  we obtain the formula.\\ ${\ }_{}^{}$ \hfill \QED

\begin{example}{\em For $\lambda=(222)$ and $l=3$, we can compute \em}
\begin{eqnarray*}
  \tilde{H}_{(2222)}(x;q,i)&=& s_{611}+s_{8}-s_{5111}-s_{71}+q^4(s_{11111111}+s_{311111}-s_{2111111}-s_{41111})+
                           \\
                           & &(q^4+1)(s_{2222}+s_{332}+s_{4211}+s_{44}-s_{431}-s_{3221})\\
                           &=&(1-q^4)(1-q^8) \ p_4\circ h_2\left (\frac{x}{1-q} \right ).                                 
\end{eqnarray*}

\end{example}

\begin{remark}
  {\em At $\ t=\zeta$, a primitive $l$-th root of unity, the inverse
    of the norm of the Macdonald polynomial $P_{(r^l)}(x;q,t)$
    satisfies
    $$
    \frac{1}{\langle P_{(r^l)}(x;q,t)\ ,\ P_{(r^l)}(x;q,t)
      \rangle_{q,t}} \Bigg{\vert}_{t=\zeta} = 0\, .$$
    Consequently, we
    obtain the following specializations \em}
$$ Q_{(r^l)}(x;q,\zeta)  = 0 \quad \text{{\em and \em}} \quad  J_{(r^l)}(x;q,\zeta) = 0.$$
\end{remark}


\section{Pieri formula at roots of unity}

In order to prove the factorization formulas, we give an auxiliary
result, in Proposition \ref{Psiroot}, on the coefficients of Pieri
formula at root of unity (c.f. \cite{macd}, Chap. VI, formula (6.24 ii))
\begin{eqnarray}\label{PieriFormula}
Q^{'}_{\mu}(x; q, t)g_r^{'}(x; q, t) =
\sum_{\lambda}\psi_{\lambda/\mu}(q,t)\ Q^{'}_{\lambda}(x; q, t),
\end{eqnarray}
with
$$
\forall r \ge 0\ ,g_r^{'}\left (\frac{1-q}{1-t}x;q,t\right ) = 
\sum_{\vert \lambda \vert = n} z_\lambda(q,t)^{-1}p_\lambda\left (\frac{1-q}{1-t}x \right ) \ . 
$$
Let $\lambda$ and $\mu$ be partitions such that
$\lambda/\mu$ is an horizontal ({\it r}-)strip $\theta$.  Let
$C_{\lambda/\mu}$ (resp.  $R_{\lambda/\mu}$) be the union of columns
(resp. rows) of $\lambda$ that intersects with $\theta$, and
$D_{\lambda/\mu}=C_{\lambda/\mu}-R_{\lambda/\mu}$ the set theoretical
difference.  Then it can be seen from the definition that for each
cell $s$ of $D_{\lambda/\mu}$ (resp.
$D_{\tilde{\lambda}/\tilde{\mu}}$) there exists a unique connected
component of $\theta$ (resp. $\tilde\theta$), which lies in the same
row as $s$.  We denote the corresponding component by $\theta_s$
(resp. $\tilde\theta_s$).\\

 \noindent Suppose that $l$ and $r$ are positive integers.  Set
 $\tilde{\lambda}=\lambda\cup(r^l)$ and $\tilde{\mu}=\mu\cup(r^l)$.
 We shall consider the difference between
 $D_{\tilde{\lambda}/\tilde{\mu}}$ and $D_{\lambda/\mu}$. It can be
 seen that there exists a projection
 $$
 p=p_{\lambda/\mu} : D_{\tilde{\lambda}/\tilde{\mu}}
 \longrightarrow D_{\lambda/\mu}.
 $$
 The cardinality of the fiber of each cell $s=(i,j)\in
 D_{\lambda/\mu}$ is exactly one or two.  Let $J_s$ denote the set of
 second coordinates of the cells in $\theta_s$.  If all elements of
 $J_s$ are all strictly larger than $r$, the fiber $p^{-1}(s)$
 consists of a single element $s=(i,j)$.  If all elements of $J_s$ are
 strictly smaller than $r$, then the fiber $p^{-1}(s)$ consists of a
 single element $\tilde{s}:=(i, j+l)$.  In the case where $J_s$
 contains $r$, then the fiber $p^{-1}(s)$ consists of exactly two
 elements $s=(i,j)$ and $\tilde{s}=(i,j+l)$.  For the case where $r\in
 J_s$, we have the following lemma, which follows immediately from the
 definition of the projection $p=p_{\lambda/\mu}$.


\begin{lemma}\label{ArmLegroot} 
  Let $s=(i,j)$ be a cell of $D_{\lambda/\mu}$ and $\tilde{s}=(i,j+l)$
  be a cell of $D_{\tilde{\lambda}/\tilde{\mu}}$ such that $r\in J_s$.
  The arm length, the arm-colength, the leg length and the
  leg-colength satisfy the following properties
\begin{eqnarray}
\label{item1} a_{\tilde{\mu}}(s) = a_{\tilde{\lambda}}(\tilde{s})\ , & \quad &  l_{\tilde{\mu}}(s)-l_{\tilde{\lambda}}(\tilde{s})=l \ , \\
\label{item2} a_{\tilde{\mu}}(\tilde{s})=a_{\mu}(s)              \ , & \quad &  l_{\tilde{\mu}}(\tilde{s})=l_{\mu}(s)               \ ,\\
\label{item3} a_{\tilde{\lambda}}(s)=a_{\lambda}(s)              \ , & \quad &  l_{\tilde{\lambda}}(s)-l_{\lambda}(s)=l             \ .
\end{eqnarray}

\end{lemma}



\begin{proposition}\label{Psiroot}
Let $\lambda$ and $\mu$ be two partitions 
such that $\mu\subset\lambda$ and 
$\theta=\lambda-\mu$ an horizontal strip. 
Let $r$ and $l$ be positive integers and  
$\zeta$ a primitive root of unity. 
It follows that
$$
\psi_{\lambda\cup (r^l)/\mu \cup (r^l)}(q, \zeta) =
\psi_{\lambda/\mu}(q, \zeta).
$$
\end{proposition}

\noindent
\proof
Recall that for a cell $s$ of the partition $\nu$,
$$
\psi_{\lambda/\mu}(q, t) = \prod_{s\in D_{\lambda/\mu}}
\frac{b_{\mu}(s)}{b_{\lambda}(s)},
$$
where
$$
b_{\nu}(s) = \frac {1-q^{a_{\nu}(s)}t^{l_{\nu}(s)+1}}
{1-q^{a_{\nu}(s)+1}t^{l_{\nu}(s)}}\ .
$$
 
\noindent If $s=(i,j)\in \lambda$ satisfies the condition 
$r<j$ for all $j\in J_s$, then the fiber $p^{-1}(s)$ of the projection
$p$ is $\{s=(i,j)\}$, and we have $ a_{\mu}(s)=a_{\tilde{\mu}}(s),
a_{\lambda}(s)=a_{\tilde{\lambda}}(s) $ and $
l_{\mu}(s)+l=l_{\tilde{\mu}}(s),
l_{\lambda}(s)+l=l_{\tilde{\lambda}}(s) $.  It is clear from these
identities that $ b_{\mu}(s)/b_{\lambda}(s) =
b_{\tilde{\mu}}(s)/b_{\tilde{\lambda}}(s) $ at $t=\zeta$ in this case.
Suppose that $s$ satisfies $j<r$ for all $j\in J_s$.  In this case,
the fiber $p^{-1}(s)$ consisits of a single element $\{\tilde{s}=(i,
j+l)\}$, and we have $ a_{\mu}(s)=a_{\tilde{\mu}}(s)$ and
$a_{\lambda}(s)=a_{\tilde{\lambda}}(s)$ and
$l_{\mu}(s)=l_{\tilde{\mu}}(s)$ and
$l_{\lambda}(s)=l_{\tilde{\lambda}}(s) $.  Hence we have $
b_{\mu}(s)/b_{\lambda}(s) = b_{\tilde{\mu}}(s)/b_{\tilde{\lambda}}(s)
$.  Consider the case where $r\in J_s$. In this case, the fiber
$p^{-1}(s)$ consists of two elements $\{s, \tilde{s}\}$. Let us consider
$$
 \prod_{u\in p^{-1}(s)}
 \frac
 {b_{\tilde{\mu}}(u)}
  {b_{\tilde{\lambda}}(u)}
  =
 \frac
        {1-q^{a_{\tilde{\mu}}(s)}t^{l_{\tilde{\mu}}(s)+1}}
  {1-q^{a_{\tilde{\mu}}(s)+1}t^{l_{\tilde{\mu}}(s)}}
 \frac
        {1-q^{a_{\tilde{\lambda}}(s)+1}t^{l_{\tilde{\lambda}}(s)}}
  {1-q^{a_{\tilde{\lambda}}(s)}t^{l_{\tilde{\lambda}}(s)+1}}
 \frac
        {1-q^{a_{\tilde{\mu}}(\tilde{s})}t^{l_{\tilde{\mu}}(\tilde{s})+1}}
  {1-q^{a_{\tilde{\mu}}(\tilde{s})+1}t^{l_{\tilde{\mu}}(\tilde{s})}}
 \frac
        {1-q^{a_{\tilde{\lambda}}(\tilde{s})+1}t^{l_{\tilde{\lambda}}(\tilde{s})}}
  {1-q^{a_{\tilde{\lambda}}(\tilde{s})}t^{l_{\tilde{\lambda}}(\tilde{s})+1}}\ .
$$
By (\ref{item1}) it follows that 
$$
\left\{ \frac {1-q^{a_{\tilde{\mu}}(s)}t^{l_{\tilde{\mu}}(s)+1}}
  {1-q^{a_{\tilde{\mu}}(s)+1}t^{l_{\tilde{\mu}}(s)}} \right\}^{-1}
\Bigg{\vert}_{t=\zeta} = \frac
{1-q^{a_{\tilde{\lambda}}(\tilde{s})+1}t^{l_{\tilde{\lambda}}(\tilde{s})}}
{1-q^{a_{\tilde{\lambda}}(\tilde{s})}t^{l_{\tilde{\lambda}}(\tilde{s})+1}}
\Bigg{\vert}_{t=\zeta}.
$$
It also follows by (\ref{item2})  
$$
        \frac
        {1-q^{a_{\tilde{\mu}}(\tilde{s})}t^{l_{\tilde{\mu}}(\tilde{s})+1}}
  {1-q^{a_{\tilde{\mu}}(\tilde{s})+1}t^{l_{\tilde{\mu}}(\tilde{s})}}\Bigg{\vert}_{t=\zeta}
 =
 \frac
        {1-q^{a_{\mu}(s)}t^{l_{\mu}(s)+1}}
  {1-q^{a_{\mu}(s)+1}t^{l_{\mu}(s)}}\Bigg{\vert}_{t=\zeta},
$$
and from (\ref{item3})
$$
        \frac
        {1-q^{a_{\tilde{\lambda}}(s)+1}t^{l_{\tilde{\lambda}}(s)}}
  {1-q^{a_{\tilde{\lambda}}(s)}t^{l_{\tilde{\lambda}}(s)+1}}\Bigg{\vert}_{t=\zeta}
 =
 \frac
        {1-q^{a_{\lambda}(s)+1}t^{l_{\lambda}(s)}}
  {1-q^{a_{\lambda}(s)}t^{l_{\lambda}(s)+1}}\Bigg{\vert}_{t=\zeta}.
$$
Therefore, it follows that 
$$
         \prod_{u\in p^{-1}(s)}
         \frac
        {b_{\tilde{\mu}}(u)}
        {b_{\tilde{\lambda}}(u)}
        =
        \frac
        {b_{\mu}(s)}{b_{\lambda}(s)}.
$$
Combining these, 
the assertion follows. 
\QED


\section{Factorization formulas}

In this section, we shall show factorization formulas for different
kinds of Macdonald polynomials at roots of unity.

\begin{theorem}\label{QPfact}
  Let $l$ be a positive integer and $\zeta$ a primitive $l$-th root of
  unity.  Let $\mu=(1^{m_1},2^{m_2},\cdots, n^{m_n})$ be a partition of a
  positive integer $n$.  For each $i$, let $m_i=lq_i+r_i$ with $0\leq
  r_i\leq l-1$ and let $\bar{\mu}=(1^{r_1}2^{r_2}\cdots m^{r_n})$.
  The function $Q_{\mu}^{'}$ satisfy the following factorisation
  formula at $t=\zeta$
        \begin{eqnarray}
                 \quad Q_{\mu}^{'}(x; q, \zeta) &=& \left(
                Q_{(1^l)}^{'}(x; q, \zeta) \right)^{q_1} \left(
                Q_{(2^l)}^{'}(x; q, \zeta) \right)^{q_2} \cdots \left(
                Q_{(n^l)}^{'}(x; q, \zeta) \right)^{q_n}
                Q_{\bar{\mu}}^{'}(x; q, \zeta). \quad
              \end{eqnarray}
\end{theorem}
        
\noindent
\proof
We shall show that the $\mathbb{C}$-linear map defined by
\begin{eqnarray*}
        f_r : \Lambda_F        & \longrightarrow & \Lambda_F \\ 
              Q_{\mu}^{'}(x; q, \zeta) & \longmapsto     & Q_{\mu\cup (r^l)}^{'}(x; q, \zeta)\ ,
\end{eqnarray*}
is an $\Lambda_{\mathbb{C}(q)}$-linear map. Let $\zeta$ be a primitive
$l$-th root of unity. From (\ref{PieriFormula}), we have
$$
Q^{'}_{\mu}(x; q, \zeta)g^{'}_k(x; q, \zeta) =
\sum_{\lambda}\psi_{\lambda/\mu}(q,\zeta)Q^{'}_{\lambda}(x; q, \zeta)\ ,
$$
where the sum is taken over the partitions $\lambda$ such that
$\lambda-\mu$ is an horizontal $k$-strip. Using the result of
Proposition \ref{Psiroot}, it follows that 
\begin{eqnarray*}
Q^{'}_{\mu\cup (r^l)}(x; q, \zeta)g^{'}_k(x; q, \zeta) &=&
\sum_{\lambda}\psi_{\lambda\cup (r^l)/\mu\cup
  (r^l)}(q,\zeta)Q^{'}_{\lambda\cup (r^l)}(x; q, \zeta) \\
& =&\sum_{\lambda}\psi_{\lambda /\mu}(q,\zeta)Q^{'}_{\lambda\cup (r^l)}(x; q, \zeta) .
\end{eqnarray*}
Consequently, for each $r\ge 1$, the multiplication by $g_k$ commutes
with the morphism $f_r$. Since the family $\lbrace g_k^{'}(x; q,
\zeta_l)\rbrace_{k\ge 1}$ generates the algebra
$\Lambda_{\mathbb{C}(q)}$ (see \cite{macd}, Chap. VI, formula (2.12)),
the map $f_r$ is $\Lambda_{{\bf C}(q)}$-linear. This implies that
\begin{eqnarray*}
 \forall F \in \Lambda_{\mathbb{C}_q},~ f_r(F(x)) & = & F(x) f_r(1)   \\
                                  & = & F(x)Q_{r^l}^{'}(x; q, \zeta).
\end{eqnarray*}
\QED

\begin{corollary}
  With the same notation as in Theorem \ref{QProot}, we have the
  following factorisation formula for the Macdonald polynomials
  $\tilde{H}_{\mu}(x;q,t)$
\begin{eqnarray}
  \tilde{H}_{\mu}(x; q, \zeta) & =& \left( \tilde{H}_{(1^l)}(x; q,
                       \zeta) \right)^{q_1} \left( \tilde{H}_{(2^l)}(x; q,
                       \zeta) \right)^{q_2} \cdots \left( \tilde{H}_{(n^l)}(x;
                       q, \zeta) \right)^{q_n} \tilde{H}_{\bar{\mu}}(x; q,
                       \zeta). \quad
\end{eqnarray}
\end{corollary}


\noindent
\proof
 If we define 
$$
        \Psi_{\lambda/\mu}(q, t)
 :=
 \psi_{\lambda/\mu}(q,t)
 \frac
        {c'_{\mu}(q, t)}{c'_{\lambda}(q, t)},
        $$
        then the Pieri formula for the modified integral form
        $J_\mu\left (\frac{x}{1-t};q,t\right )$ is written as follows
$$
    J_\mu\left (\frac{x}{1-t};q,t\right )g^{'}_r\left ( \frac{x}{1-q}\right )    
 =
 \sum_{\lambda}
 \Psi_{\lambda/\mu}(q,t) J_{\lambda}\left (\frac{x}{1-t}; q, t\right ),
$$
where the sum is over the partitions $\lambda$ 
such that $\lambda-\mu$ is a  horizontal $k$-strip. \\
Let a positive integer $r$ be arbitrarily fixed, 
and $\tilde{\nu}$ denote the partition $\nu\cup(r^l)$. 
Since we have already shown that 
$
        \psi_{\tilde{\lambda}/\tilde{\mu}}(q, \zeta)
 =
 \psi_{\lambda/\mu}(q, \zeta)
$, 
it suffices to show that
$$
        \frac
 {      c'_{\tilde{\mu}}(q, \zeta)      }
 {      c'_{\tilde{\lambda}}(q, \zeta)  }
 =
 \frac
 {      c'_{\mu}        (q, \zeta)              }
 {      c'_{\lambda}    (q, \zeta)      }.
$$
We shall actually show that
$$
        \frac
 {      c'_{\tilde{\mu}}(q, \zeta)      }
 {      c'_{\mu}        (q, \zeta)              }
 =
 \frac
 {      c'_{\tilde{\lambda}}(q, \zeta)  }
 {      c'_{\lambda}    (q, \zeta)      }.
$$
        
\noindent It follows from the definition that
\begin{eqnarray*}
        \frac
        {c'_{\tilde{\mu}}(q,t)}
  {c'_{\mu}(q,t)}
        &=&
        \frac
        {\prod_{s\in\tilde{\mu}}(1-q^{a_{\tilde{\mu}}(s)+1}t^{l_{\tilde{\mu}}(s)})}
     {\prod_{s\in\mu}(1-q^{a_{\tilde{\mu}}(s)+1}t^{l_{\tilde{\mu}}(s)})}\\
   &=&
        \frac
        {
        \prod_{s\in\tilde{\mu}\atop s\notin(r^l)}
     (1-q^{a_{\tilde{\mu}}(s)+1}t^{l_{\tilde{\mu}}(s)})
    }
        {
        \prod_{s\in\mu}
     (1-q^{a_{\tilde{\mu}}(s)+1}t^{l_{\tilde{\mu}}(s)})
    }
    \prod_{s\in (r^l)\subset\tilde{\mu}}
    (1-q^{a_{\tilde{\mu}}(s)+1}t^{l_{\tilde{\mu}}(s)}).
\end{eqnarray*}
The Young diagram of the partition $\tilde{\mu}$ 
is the disjoint union of the cells 
$
        \{\tilde{s}\in\tilde{\mu}|s\in\mu\}
$ 
and $(r^l)$. 
For each $s\in \mu$, we have as seen in previous theorem that 
$a_{\tilde{\mu}}(\tilde{s})=a_{\mu}(s)$, 
and $l_{\tilde{\mu}}(\tilde{s})=l_{\mu}(s)$ or $l_{\mu}(s)+l$. 
Hence at $t=\zeta$, we have
\begin{eqnarray}
                \frac
        {c'_{\mu}(q,\zeta)}
  {c'_{\tilde{\mu}}(q,\zeta)}
  =
  \prod_{s\in(r^l)\subset\tilde{\mu}}
  (1-q^{a_{\tilde{\mu}}(s)+1}t^{l_{\tilde{\mu}}(s)}) \quad \text{ and }
\end{eqnarray}
\begin{eqnarray}
        \frac
        {c'_{\lambda}(q,\zeta)}
 {c'_{\tilde{\lambda}}(q,\zeta)}
 =
 \prod_{s\in(r^l)\subset\tilde{\lambda}}
 (1-q^{a_{\tilde{\lambda}}(s)+1}t^{l_{\tilde{\lambda}}(s)}).
\end{eqnarray}

Although there is a difference between the positions where the block $(r^l)$ is inserted 
in the Young diagram of $\mu$ and $\lambda$, 
(3.1) and (3.2) coincide at $t=\zeta$, 
since $a_{\tilde{\mu}}(s)=a_{\tilde{\lambda}}(s)$ for each $s\in (r^l)$. 
Thus we have
$$
        \frac{c'_{\mu}(q, \zeta)}{c'_{\tilde{\mu}}(q, \zeta)}
 =
 \frac{c'_{\lambda}(q, \zeta)}{c'_{\tilde{\lambda}}(q, \zeta)}.
$$
\QED

\noindent Let $\nu=(\nu_1,\ldots,\nu_p)$ be a partition. For some $l\ge 0$, we
denote by $\nu^l$ the partition where each part of $\nu$ is repeated
$l$ times. We can give a more explicit expression for the
factorisation formula in the special case where $\mu=\nu^l$.

\begin{corollary}
  Let $\nu$ be a partition and $l$ a positive integer. We have the
  following special cases for the factorisation formulas
\begin{eqnarray}
 Q^{'}_{\nu^l}(X;q,\zeta)     & = & (-1)^{(l-1)\vert \nu \vert}p_l\circ h_\nu(x)\ ,\\ 
 \tilde{H}_{\nu^l}(X;q,\zeta) & = & \prod_{j=1}^{l(\nu)} \prod_{i=1}^{\nu_j}(q^{il}-1)\ p_l\circ h_\nu\left(\frac{x}{1-q}\right )\ .
\end{eqnarray}
\end{corollary}

\begin{example} {\em For $\lambda =(222111)$ and $k=3$, we can compute the specialization \em}
  \begin{eqnarray*}
    Q^{'}_{222111}(x;q,e^{2i\pi/3})&=&-s_{22221}-s_{321111}+s_{3222} +
   s_{33111}-s_{3321}+ 3 s_{333}+ s_{411111} \\
& & - 2 s_{432} + 2 s_{441} - 
   s_{51111} + 2 s_{522} - 2 s_{54} + s_{6111} - 2 s_{621} + 2 s_{63}\\
& & + s_{711} -  s_{81} + s_{9} + s_{222111}\\
&=&p_{3}\circ h_{21}(x). \\
  \end{eqnarray*}
\end{example}


\section{A generalization of the plethystic formula}

In this section, using the factorisation formula given in Theorem
\ref{QPfact}, we shall give a generalization of the plethystic formula
obtained by specializing Macdonald polynomials at roots of unity in
Theorem \ref{QProot}. For $\lambda$ a partition, let consider the
following map which is the plethystic substitution by the powersum
$p_\lambda$
\begin{eqnarray*}
{\Psi}_\lambda: {\Lambda}_F & \longrightarrow & {\Lambda}_F            \\
                 f          & \longmapsto & p_\lambda\circ f \ .
\end{eqnarray*}

\begin{lemma}\label{plethysm}
  Let $\lambda$ and $\mu$ be two partitions, the maps ${\Psi}$ satisfy
  the multiplicative rule
  $$
  {\Psi}_\lambda \left ( f \right ) {\Psi}_\mu \left ( f \right) =
  {\Psi}_{\lambda\cup \mu} \left ( f \right )\ . $$
\end{lemma}


\begin{proposition} 
Let $d$ be an integer such that $d\vert l$ and $\zeta_d$ be a
primitive $d$-th root of unity,
\begin{eqnarray}
 Q^{'}_{(r^l)}(x;q,\zeta_d)&=&(-1)^{\frac{rl(d-1)}{d}}
p_d^{l/d}\circ h_r(x)\ .
\end{eqnarray}
\end{proposition}

\noindent
\proof Let $d$ and $l$ be two integers such that $d$ divide
$l$. Let $\mu=(r^l)$ the rectangle partition with parts of length $r$.
Using the factorisation formula described in Theorem \ref{QPfact}, we can write
\begin{eqnarray}\label{QPRootd}
Q^{'}_{(r^l)}(x;q,\zeta_d)=\left ( Q^{'}_{(r^d)}(x;q,\zeta_d) \right
)^{l/d}.
\end{eqnarray}
With the specialization formula at root of
unity written in Theorem \ref{QProot}, we have

\begin{eqnarray*}
\left ( Q^{'}_{(r^d)}(x;q,\zeta_d) \right)^{l/d}
&=&
\left (
(-1)^{(d-1)r}p_d\circ h_r(x)\right )^{l/d} \\
&=&
(-1)^{\frac{lr(d-1)}{d}}\left (p_d\circ h_r(x)\right )^{l/d}
\end{eqnarray*}
Using the Lemma \ref{plethysm}, we obtain
$$\left ( Q^{'}_{(r^d)}(x;q,\zeta_d)
\right)^{l/d}=(-1)^{\frac{lr(d-1)}{d}}p_d^{l/d}\circ h_r(x)\ . $$
Finally, we obtain by the factorization formula of Theorem 5.1
$$ Q^{'}_{(r^l)}(x;q,\zeta_d)=(-1)^{\frac{lr(d-1)}{d}}p_d^{l/d}\circ h_r(x)\ . $$
\QED

\noindent Using the same proof, we can write a similar specialization for
integral forms of the Macdonald Polynomials.
\begin{corollary}
  With the same notations as in Proposition 6.1, the modified
  Macdonald polynomials $\tilde{H}_\lambda(x;q,t)$ satisfy the same
  specialization
\begin{eqnarray}
\tilde{H}_{(r^l)}(x;q,\zeta_d)&=& \prod_{i=1}^{r}(q^{il}-1)
p_d^{l/d}\circ h_r\left (\frac{x}{1-q}\right )\ .
\end{eqnarray}
\end{corollary}

\begin{example}{\em For $\lambda=(222222)$, i.e $r=2$ and $l=6$ and $d=3$, we can compute \em}
  \begin{eqnarray*}
  {Q}^{'}_{(222222)}(x;q,e^{2i\pi/3} ) &=& 
 - s_{3 2 2 2 2 1} + s_{3 3 2 2 2} + 2 s_{3 3 3 1 1 1} -
   2 s_{3 3 3 2 1} + 2 s_{3 3 3 3} + s_{4 2 2 2 1 1} -
   2 s_{4 3 2 1 1 1} \\ & & + s_{4 3 2 2 1} + 2 s_{4 4 1 1 1 1} -
   s_{4 4 2 2} + 4 s_{4 4 4} + s_{5 2 2 1 1 1} -
   2 s_{5 2 2 2 1} + s_{5 3 2 1 1} - 2 s_{5 4 1 1 1} \\ & &+
   s_{5 4 2 1} - 4 s_{5 4 3} + 3 s_{5 5 2} - s_{6 2 1 1 1 1} +
   2 s_{6 2 2 2} + s_{6 3 1 1 1} - 2 s_{6 3 2 1} + 4 s_{6 3 3} \\ & &+
   s_{6 4 1 1} - 3 s_{6 5 1} + 3 s_{6 6} + s_{7 1 1 1 1 1} -
   2 s_{7 3 2} + 2 s_{7 4 1} - s_{8 1 1 1 1} + 2 s_{8 2 2} -
   2 s_{8 4}\\ & & + s_{9 1 1 1} - 2 s_{9 2 1} + 2 s_{9 3} + s_{10 1 1} -
   s_{11 1} + s_{12} + s_{2 2 2 2 2 2}\\
&=& p_3^2\circ h_2(x)=p_{(33)}\circ h_2(x).
\end{eqnarray*}
\end{example}


\section{Macdonald polynomials at roots of unity and cyclic 
characters of the symmetric group}

In the following, we will denote the symmetric group of order $k$ by
$\mathfrak{S}_k$. Let $\Gamma \subset \mathfrak{S}_k$ be a cyclic
subgroup generated by an element of order $r$. As $\Gamma$ is a
commutative subgroup its irreducible representations are
one-dimensional vector spaces. The corresponding maps $(\gamma
_j)_{j=0\ldots r-1}$ can be defined by
\begin{eqnarray*}
\gamma_j : & \Gamma  \longrightarrow & GL(\mathbb{C})\ \simeq\ \mathbb{C}^{*} \\
           & \tau    \longmapsto     & \zeta_r^j\ ,
\end{eqnarray*}
where $\zeta_r$ is a $r$-th primitive root of unity (See \cite{serre}
for more details). In \cite{Foulkes}, Foulkes considered the Frobenius
characteristic of the representations of $\mathfrak{S}_k$ induced by
these irreducible representations and obtained an explicit formula
that we will give in the next proposition. 
\\Let $k$ and $n$ be two
positive integers such that $u=(k,d)$ (the greater common divisor
between $k$ and $n$) and $d=u\cdot m$. Let us define the Ramanujan (or Von
Sterneck) sum $c(k,d)$ by
$$c(k,d)=\frac{\mu(m)\phi(d)}{\phi(m)}$$
where $\mu$ is the Moebius
function and $\phi$ the Euler totient. The quantity $c(k,d)$
corresponds to the sum of the $k$-th powers of the primitive $d$-th
roots of unity (the previous expression was given first by H\"older in
\cite{HardyWright}).

\begin{proposition}
  Let $\tau$ be a cyclic permutation of length $k$ and $\Gamma$ the
  maximal cyclic subgroup of $\mathfrak{S}_k$ generated by $\tau$. Let
  $j$ be a positive integer less than $k$. The Frobenius
  characteristic of the representation of $\mathfrak{S}_k$ induced by
  the irreducible representation of $\Gamma$, $\gamma_j: \tau
  \longmapsto \zeta_r^j$, is given by
\begin{eqnarray}
l_k^{(j)}(x)=\frac{1}{k}\sum_{d \vert k} c(j,d)~p_d^{k/d}(x).
\end{eqnarray}
\end{proposition}

\begin{example} 
  For $\mathfrak{S}_6$ and $k=2$, the cyclic character $l^{(2)}_6$
  expanded on powersums and Schur basis is
  \begin{eqnarray*}
  l_6^{(2)}&=&\frac{1}{6}\left ( p_{111111}+p_{222}-p_{33}-p_{6}
  \right ) \\
  &=& s_{51}+2s_{42}+s_{411}+3s_{321}+2s_{3111}+s_{222}+s_{2211}+s_{21111}.
  \end{eqnarray*}
\end{example}

\begin{theorem}
  Let $r$ and $l$ be two positive integers. The specialization of the
  Macdonald polynomials indexed by the rectangle partition $(r^l)$ at
  a primitive $l$-th root of unity is equivalent to
\begin{eqnarray}\label{QPMod}
  Q^{'}_{(r^l)}(x;q,t)\mod \ \Phi_l(t)\ &=&\ \sum_{j=0}^{l-1}t^j~(l_l^{(j)}\circ
  h_r)(x)\ .
\end{eqnarray}
\end{theorem}
\noindent
\proof We will first give a generalization of the Moebius
inversion formula due to E. Cohen (see \cite{Cohen} for the original
work and \cite{Desar} for a simpler proof). Let
$$P(q)=\sum_{k=0}^{n-1}a_kq^k \ ,$$
be a polynomial of degree less than
$n-1$ with coefficients $a_k$ in $\mathbb{Z}$. $P$ is said to be even
modulo $n$ if $$(i,n)=(j,n)\ \Longrightarrow \ a_i=a_j.$$
\begin{lemma}\label{cohen} The polynomial $P$ is even modulo $n$ if and only if for every 
  divisor $d$ of $n$, the residue of $P$ modulo the $d$-th cyclotomic
  polynomial $\Phi_d$ is a constant $r_d$ in $\mathbb{Z}$. In this
  case, one has
  $$
  a_k = \frac{1}{n}\sum_{d\vert n}c(k,d)\ r_d \quad \text{ and }
  \quad r_d = \sum_{t\vert n}c(n/d,t)\ a_{n/t}.
$$
\end{lemma}

\noindent Let $d$ be an integer such that $d\vert l$. By expanding $Q^{'}_{(r^l)}(x;q,t)$
(and more generally $Q^{'}_\lambda(x;q,t)$) on the Schur basis, we can
define a kind of $(q,t)$-Kostka polynomials $K^{'}_{\mu, (r^l)}(q,t)$
by
$$Q^{'}_{(r^l)}(x;q,t)=\sum_{\mu} K^{'}_{\mu, (r^l)}(q,t)~ s_\mu(x).$$

\noindent Let $\mu$ be a partition and $d$ an integer such that $d\vert l$. The polynomial 
$P_{\mu}^q(t)=\sum_{j=0}^{l-1}a_j(q)t^j$ is the residue modulo
$1-t^l$ of the $(q,t)$-Kostka polynomial $K^{'}_{\mu,(r^l)}(q,t)$, if
and only if, for all $\zeta_d$ primitive $d$-th root of unity,
$$P_\mu^q(\zeta_d)=K^{'}_{\mu,(r^l)}(q,\zeta_d).$$  Using
Theorem \ref{QPfact} , one has
$$
P_\mu^q(\zeta_d)=(-1)^{(d-1)rl/d}\ \langle \ p_d^{l/d}\circ h_r(x)\ 
,\ s_\mu(x)\ \rangle.$$
Consequently $P(\zeta_d)$ is an integer since the entries of the 
transition matrix between the powersums and the Schur functions are all
integers. Using the Lemma \ref{cohen}, we obtain
\begin{eqnarray*}
a_j(q)&=&\frac{1}{l}\sum_{d \vert l}c(j,d)\ \langle \ p_d^{l/d}\circ h_r(x)\ ,\ 
s_\mu(x)\ \rangle \\
   &=&\langle \ l^{(j)}_{l}\circ h_r(x)\ ,\ s_\mu(x) \ \rangle.
\end{eqnarray*}
\QED

\begin{example}
Let define $g$ as the right hand side of \ref{QPMod} for $r=2$ and $l=3$
 \begin{eqnarray*}
g(t) &=& l_3^{(0)}\circ h_2 + t\ l_3^{(1)}\circ h_2 + t^2\ l_3^{(2)}\circ h_2 \\
  &=& s_{411} + (t^2 + t + 1) s_{42} + t(t + 1)s_{51} + t(t + 1)s_{321} + s_{222} + s_{33} + s_6\ .
 \end{eqnarray*}
 The specialization of $g(t)$ and $Q^{'}_{222}(X;q,t)$ at $t$ the 3-th
 primitive roots of unity satisfy
\begin{eqnarray*}
 g(j)   = Q^{'}_{222}(X;q,j)   & = & s_{411}-s_{51}-s_{321}+s_{222}+s_{33} + s_6\ , \\
 g(j^2) = Q^{'}_{222}(X;q,j^2) & = & s_{411}-s_{51}-s_{321}+s_{222}+s_{33} + s_6\ .
\end{eqnarray*}
 
\end{example}

\begin{corollary}
  For two positive integers $r$ and $l$, the same residue formulas
  occurs for the modified Macdonald polynomials
  $\tilde{H}_{(r^l)}(x;q,t)$
  \begin{eqnarray}
\tilde{H}_{(r^l)}(x;q,t)\mod \Phi_l(t) \ &\equiv & \prod_{i=1}^{r}(q^{il}-1)\ 
  \sum_{j=0}^{l-1}t^j~(l_l^{(j)}\circ h_r)\left (\frac{x}{1-q}\right )\ .
\end{eqnarray}
\end{corollary}

\section{Congruences for $(q,t)$-Kostka polynomials}
For a given partition partition $\lambda$, let denote by
$\tilde{s}^{(q)}_\lambda$ the symmetic function defined by
$$
\tilde{s}^{(q)}_\lambda(x)=s_\lambda\left ( \frac{x}{1-q} \right
).$$
Let $\star$ be the internal product on $\Lambda_{F}$ defined by
(see \cite{macd}, Chap. I, formula (7.12))
$$ p_\lambda \star p_\mu = \delta_{\lambda,\mu}z_\lambda p_\lambda .$$

\begin{proposition}
Let $r$ and $l$ be two positive integers and $\mu$ a partition of
weight $nl$. Let denote by $\Phi_l(t)$ the cyclotomic polynomial of
order $l$. The $(q,t)$-Kostka polynomial $\tilde{K}_{(r^l),\mu}(q,t)$
satisfy the following congruence modulo $\Phi_l(t)$
\begin{eqnarray}
 \tilde{K}_{\mu, (r^l)}(q,t) \equiv \prod_{i=1}^{r}(q^{il}-1)\
\widetilde{s}^{\ (q)}_\mu(1,t,t^2,\ldots,t^{l-1}) \mod \Phi_l(t).
\end{eqnarray}
More generally, for all partitions $\nu=(\nu_1,\ldots,\nu_p)$ of weight $r$\ ,
\begin{eqnarray}
\tilde{K}_{\mu, \nu^l}(q,t)\equiv \prod_{j=1}^{l(\nu)}\prod_{i=1}^{\nu_j}(q^{il}-1)\
\widetilde{h_{l\nu}\star s_\mu}^{(q)}(1,t,t^2,\ldots,t^{l-1}) \mod \Phi_l(t)\ ,
\end{eqnarray}
where $l\nu=(l\nu_1,\ldots,l\nu_p)$.
\end{proposition}

\proof
  Let $\zeta$ be a primitive root of unity and
  $Z_l=\{1,\zeta,\ldots,\zeta^{l-1}\}$ be the alphabet of the
  $l$-roots of unity. Using $\lambda$-ring notations (see \cite{L} for
  more details) and Theorem \ref{QProot}, we have for all positive integer $r$
  $$
  \tilde{H}_{r^l}(X;q,\zeta)=\prod_{i=1}^{r}(q^{il}-1)(p_l\circ
  h_r)\left (\frac{x}{1-q}\right )=h_{lr}\left (\frac{Z.x}{1-q}\right )\ .$$
      Consequently, for
      all partitions $\mu$ of size $rl$, we can write
\begin{eqnarray*}
\tilde{K}_{\mu, r^l}(q,\zeta)&=& \prod_{i=1}^{r}(q^{il}-1)\
\widetilde{s}^{\ (q)}_{\mu}(Z_l)\ ,
\end{eqnarray*} 
which is equivalent to the first statement of the theorem. The second
statement follows from the following identity
$$
(p_l \circ h_{\nu})\left( \frac{x}{1-q}\right ) = (h_{l\nu} \star
h_{lr}) \left ( \frac{Z_l x}{1-q} \right )\ .$$

\QED

\begin{example}
Let consider $r=2$ and $l=3$. For $\mu=(222)$, we have 
$$ \tilde{K}_{222,222}(q, j) =  \tilde{K}_{222,222}(q, j^2) = 1+q^3 \quad \text{and} \quad 
(1-q^6)(1-q^3)\widetilde{s}^{\ (q)}_{222}(1,j,j^2) = 1+q^3\ .$$
\end{example}

\end{document}